\documentclass[preprint]{siamart1116}

\usepackage{framed}
\usepackage{enumitem}
\usepackage{amssymb}

\usepackage{braket,amsfonts}

\usepackage{array}

\usepackage[caption=false]{subfig}
\usepackage{pgfplots}

\newsiamthm{claim}{Claim}
\newsiamremark{rem}{Remark}
\newsiamremark{expl}{Example}
\newsiamremark{hypothesis}{Hypothesis}
\crefname{hypothesis}{Hypothesis}{Hypotheses}

\usepackage{algorithmic}

\usepackage{graphicx,epstopdf}

\Crefname{ALC@unique}{Line}{Lines}

\usepackage{amsopn}

\numberwithin{theorem}{section}

\usepackage{xspace}
\usepackage{bold-extra}
\usepackage[most]{tcolorbox}

\colorlet{texcscolor}{blue!50!black}
\colorlet{texemcolor}{red!70!black}
\colorlet{texpreamble}{red!70!black}
\colorlet{codebackground}{black!25!white!25}

\lstdefinestyle{siamlatex}{%
  style=tcblatex,
  texcsstyle=*\color{texcscolor},
  texcsstyle=[2]\color{texemcolor},
  keywordstyle=[2]\color{texemcolor},
  moretexcs={cref,Cref,maketitle,mathcal,text,headers,email,url},
}

\tcbset{%
  colframe=black!75!white!75,
  coltitle=white,
  colback=codebackground, 
  colbacklower=white, 
  fonttitle=\bfseries,
  arc=0pt,outer arc=0pt,
  top=1pt,bottom=1pt,left=1mm,right=1mm,middle=1mm,boxsep=1mm,
  leftrule=0.3mm,rightrule=0.3mm,toprule=0.3mm,bottomrule=0.3mm,
  listing options={style=siamlatex}
}

\newtcblisting[use counter=example]{example}[2][]{%
  title={Example~\thetcbcounter: #2},#1}

\newtcbinputlisting[use counter=example]{\examplefile}[3][]{%
  title={Example~\thetcbcounter: #2},listing file={#3},#1}

\DeclareTotalTCBox{\code}{ v O{} }
{ 
  fontupper=\ttfamily\color{black},
  nobeforeafter,
  tcbox raise base,
  colback=codebackground,colframe=white,
  top=0pt,bottom=0pt,left=0mm,right=0mm,
  leftrule=0pt,rightrule=0pt,toprule=0mm,bottomrule=0mm,
  boxsep=0.5mm,
  #2}{#1}

\patchcmd\newpage{\vfil}{}{}{}
\newcommand{\ep}[1]{{\color{black} #1}}

\definecolor{bggreen}{rgb}{0.0, 0.75, 0.0}
\newcommand{\bg}[1]{{\color{black} #1}}

\definecolor{darkviolet}{rgb}{0.58, 0.0, 0.83}
\newcommand{\bgcom}[1]{} 

\newcommand{\com}[1]{} 

\begin{tcbverbatimwrite}{tmp_\jobname_header.tex}
  \title{An improved Krylov eigenvalue strategy using the FEAST algorithm with inexact system solves}

\author{Brendan Gavin%
  \thanks{University of Massachusetts Amherst Department of Electrical and Computer Engineering
    (\email{bgavin@ecs.umass.edu}, \email{polizzi@ecs.umass.edu})}%
  \and
  Eric Polizzi\footnotemark[1]
}

\headers{IFEAST: an outer-inner iterative strategy for the FEAST eigensolver}
{Brendan Gavin and Eric Polizzi}
\end{tcbverbatimwrite}
\input{tmp_\jobname_header.tex}

\ifpdf
\hypersetup{ pdftitle={IFEAST}} 
\fi

\begin{document}
\maketitle

\begin{tcbverbatimwrite}{tmp_\jobname_abstract.tex}
\begin{abstract}
  The FEAST eigenvalue algorithm is a subspace iteration algorithm that uses contour integration in the complex plane to obtain the eigenvectors of a matrix for the eigenvalues that are located in any user-defined search interval. By computing small numbers of eigenvalues in specific regions of the complex plane, FEAST is able to naturally parallelize the solution of eigenvalue problems by solving for multiple eigenpairs simultaneously. The traditional FEAST algorithm is implemented by directly solving collections of shifted linear systems of equations; in this paper, we describe a variation of the FEAST algorithm that uses iterative Krylov subspace algorithms for solving the shifted linear systems inexactly. We show that this iterative FEAST algorithm (which we call IFEAST) is mathematically equivalent to a block Krylov subspace method for solving eigenvalue problems. By using Krylov subspaces indirectly through solving shifted linear systems, rather than directly for projecting the eigenvalue problem, IFEAST is able to solve eigenvalue problems using very large dimension Krylov subspaces, without ever having to store a basis for those subspaces. IFEAST thus combines the flexibility and power of Krylov methods, requiring only matrix-vector multiplication for solving eigenvalue problems, with the natural parallelism of the traditional FEAST algorithm. We discuss the relationship between IFEAST and more traditional Krylov methods, and provide numerical examples illustrating its behavior. 
\end{abstract}

\begin{keywords}
FEAST, contour integration, eigenvalue problem, Kryov, Arnoldi, linear system
\end{keywords}

\begin{AMS}
  65F15,65F10,15A18
\end{AMS}
\end{tcbverbatimwrite}
\input{tmp_\jobname_abstract.tex}

\section{Introduction}

Eigenvalue problems are a staple of basic linear algebra \cite{golub2012matrix,saad1992numerical}, and they underlie a wide variety of practical computing techniques. Of particular interest are problems such as ground state quantum chemistry, linear time-dependent systems, and dimensionality reduction for data sets.

Conventional algorithms for small dimension eigenvalue problems (such as QR iterations) are generally \ep{unable to cope with} the computational demands of the larger problem sizes \ep{found in}  modern applications. Iterative algorithms that are designed specifically
for approximating the solutions to large eigenvalue problems, such as Krylov subspace methods (e.g. Lanczos and Arnoldi), tend to fare much better, and these are the primary methods that are used in solving the largest eigenvalue problems in contemporary research. These methods, however, are not necessarily the most appropriate ones \ep{for modern computing architectures}, particularly as scientific computing continues to approach the exascale. Modern high performance computing architectures achieve their promise of high performance through immense parallelism; Krylov subspace methods, on the other hand, are inherently serial algorithms that happen to be able to benefit from having large amounts of memory available. Although they can be implemented and run on parallel computers, they are not able to take full advantage of parallelism by actually dividing the task at hand into a collection of smaller, independent problems. 

The likely best way forward for solving eigenvalue problems on modern parallel computing architectures
is to \bg{use} spectral slicing along with filtering techniques 
\bg{in order to divide the spectrum of a matrix} into
an arbitrary number of smaller, non-intersecting regions in the complex plane.
The eigenvalues (and corresponding eigenvectors) in each region can be filtered from the original problem and then
solved independently of those in the other regions.
As a result, one can solve for a large number of eigenvalue/eigenvector pairs in a genuinely parallel fashion.

In this paper we discuss a modification of the FEAST algorithm (which is an example of a spectral filtering technique)
that allows one to solve eigenvalue problems for large numbers of eigenvalue/eigenvector pairs by using only matrix-vector
multiplication, in order to provide a robust and naturally parallel alternative to traditional Krylov \ep{iterative}
methods \ep{for the eigenvalue problem}.
We show that this modified FEAST algorithm, which we call Iterative FEAST (IFEAST), converges linearly to the desired
eigenpairs anywhere in the spectrum, and that it is mathematically equivalent to a restarted Krylov subspace method.
Unlike other restarted Krylov subspace eigenvalue algorithms, however, IFEAST provides a clear condition for convergence when restarting, and it 
can be implemented without having to store a basis for the Krylov subspace.


\subsection{The FEAST Algorithm}
\bg{FEAST \cite{polizziFEAST,tangFEAST} uses a spectral filtering technique that can select
the eigenpairs of interest by using an approximate spectral projector combined with a subspace iteration procedure.}
It can be used to solve the generalized eigenvalue problem

\begin{equation}
\label{eqn:eigenvalueproblem}
AX_I=BX_I\Lambda_I,
\end{equation}
with
\begin{equation}
A ,B\in \mathbb{C}^{n \times n},\ X_I=\{x_1,\dots,x_m\}_{n\times m},\ \Lambda_I=diag(\lambda_1,\dotsm\lambda_m), 
\end{equation}

\noindent \ep{by finding all the $m$ eigenvectors $x_i$ whose eigenvalues $\lambda_i$}  lie in some user-defined
region in the complex plane.
For the sake of simplicity we consider only regions that are intervals ${\cal I}=(\lambda _{min},\lambda _{max} )$
on the real number line,
thereby restricting our attention to Hermitian matrices. In general, though, FEAST, and all of the results in this paper,
can be extended to
non-Hermitian matrices as well \cite{tang2014new,kestyn2016feast}. We also restrict our attention
primarily
to the standard eigenvalue problem case (i.e. $B=I$); the reasons for this will be addressed in Section \ref{sec:useGMRES}.

FEAST selects the eigenvalues to solve for by using an approximation for the spectral projector $\rho(A)=X_IX_I^T$ in order to
form a subspace $Q$ from a (possibly random) initial guess for the eigenvectors $X$, thus guaranteeing that the columns
of $Q$ span
only the eigenvectors of interest. Because $X_I$ is unknown before solving the eigenvalue problem, FEAST uses complex contour
integration in order to form an operator that is equal to $\rho(A)$:

\begin{align}
\label{eqn:contourprojector}
Q &= \rho(A)X=(X_IX_I^T)X=\frac{1}{2\pi i}\oint_{\cal C} (zI-A)^{-1}Xdz.
\end{align}

This integral can not be evaluated exactly; in practice, multiplication by $\rho(A)$ is evaluated
approximately by using a quadrature rule:
\begin{align}
\rho(A)X&=\frac{1}{2\pi i}\oint_{\cal C} (zI-A)^{-1}Xdz \label{eqn:contour_int} \\
&\approx \sum_{k=1}^{n_c} \omega _k (z_kI-A)^{-1}X\equiv\hat \rho (A) X\label{eqn:contour_int_quad}
\end{align}
\noindent The spectral projector is thus applied in an approximate way to the estimated subspace \bg{$X$} by solving
$n_c$ shifted linear systems,
and adding their solutions together in a weighted sum. Thereafter,
the original eigenvalue problem is solved approximately in the subspace spanned by $Q$ by using the Rayleigh-Ritz
procedure,
giving new estimates of the desired eigenvalues and and eigenvectors.
The estimated eigenvectors and eigenvalues are improved iteratively by repeating this procedure until convergence.

FEAST can be interpreted as a subspace iteration that uses  the approximate spectral projection operator
$\hat \rho (A)$ as a rational
  filtering/selection function:
\begin{equation}\label{eqn:rhoh}
\hat \rho (A)=\sum_{k=1}^{n_c} \omega _k (z_kI-A)^{-1}=X\hat\rho(\Lambda)X^H,
\end{equation}
where $\hat\rho(\Lambda)$ acts on each eigenvalue individually, i.e. $\hat\rho(\lambda_j)=\sum_{k=1}^{n_c} \omega _k (z_k-\lambda_j)^{-1}$. \bg{At the limit of large $n_c$, $\hat\rho(\lambda_j)$ is either equal 1 if $\lambda_j$ is inside $\cal C$, or is equal to 0 if $\lambda_j$ is outside $\cal C$.}

Like conventional subspace iterations, the convergence of FEAST is linear \cite{tangFEAST}. It is similar to shift-invert subspace iterations but, unlike a traditional shift-invert subspace iteration algorithm, FEAST uses multiple shifts to accelerate convergence, the weights and locations of which are determined in an optimal way by using complex contour integrations. The rate of convergence is both related to the size of the search subspace and to the accuracy with which the original integral in equation (\ref{eqn:contour_int}) is approximated; the more linear systems that we solve for the quadrature rule (\ref{eqn:contour_int_quad}), the better the integral is approximated, and the fewer FEAST subspace iterations are required to converge to the desired level of accuracy.
One of the benefits of this is that, because the linear systems can be solved independently of each other,
the use of additional parallel
processing power can be translated directly into a faster convergence rate simply by solving more linear systems
in parallel.
\ep{Algorithm \ref{algo:feast} summarizes the basic FEAST procedure for solving the standard Hermitian interior
  eigenvalue problem.}

\begin{algorithm}
\ep{
  \medskip
\noindent \begin{minipage}{\linewidth}
{\footnotesize

\noindent {\bf Start with:}
\begin{itemize}
\item Matrix $A \in \mathbb{C}^{n \times n}$
\item Interval ${\cal I}=(\lambda _{min},\lambda _{max} )$ wherein
  fewer than $m_0$ eigenvalues are expected to be found, and closed contour $\cal{C}$ that encloses
  $\cal I$ in the complex plane
\item Initial guess $X^{(0)} \in \mathbb{C}^{n\times m_0}$ for the search subspace spanned by the solution
  to the eigenvalue problem
\item Set of $n_c$ quadrature weights and points $(\omega_k,z_k)$ for numerically integrating equation
  (\ref{eqn:contour_int}) \footnote{Any quadrature rule can be used, e.g. Gaussian quadrature, trapezoidal
  or Zolotarev rule \cite{zolotarev}. For an explicit example of how to integrate (\ref{eqn:contour_int}) numerically,
  see \cite{polizziFEAST}.}
\end{itemize}
\medskip
\noindent {\bf For each} subspace iteration $i$:
\medskip
\begin{description}
\setlength{\itemsep}{3pt}
\item[1.] Directly solve $n_c$ shifted linear systems for $Y_k^{(i)}\in \mathbb{C}^{n\times m_0}$.
\begin{equation}\nonumber
(z_kI-A)Y_k^{(i)}=X^{(i)}, \ \ 1\leq k\leq n_c
\end{equation}
\item[2.] Form the filtered subspace $Q$
\begin{equation}\nonumber
Q=\hat \rho(A)X^{(i)} = \sum_{k=1}^{n_c}\omega_k Y_k^{(i)}
\end{equation}

\item[3.] Perform Rayleigh-Ritz procedure to find a new estimate for eigenvalues and eigenvectors:
\begin{description}
\setlength{\itemsep}{3pt}
\item[i.] Solve the generalized reduced eigenvalue problem for $X_Q  \in \mathbb{C}^{m_0\times m_0} $
\begin{equation} \nonumber
  A_Q X_Q= B_Q X_Q \Lambda
\end{equation}
\noindent with $A_Q =Q^T A Q \text{ and } B_Q=Q^TQ$
\item[ii.] Get new estimate for subspace $X^{(i+1)}$: $X^{(i+1)} = Q X_Q$
\end{description}
\item[4.] Calculate the FEAST eigenvector residual $||R_F||=max\  ||Ax_j-\lambda _j x_j||,\ 1 \leq j
  \leq m_0,\ \lambda_j \in {\cal I}$. If $||R_F||$ is above a given tolerance, {\bf GOTO 1}.
\end{description}
}
\end{minipage}
\caption{\label{algo:feast} The FEAST Hermitian algorithm for solving $AX_I=X_I\Lambda_I$}
}
\end{algorithm}

\subsection{Challenges for FEAST}

FEAST is most useful when applied to sparse matrices of high dimension. In this case, one would
typically use an optimized sparse
direct solver (such as PARDISO \cite{pardisosite}) for the solution of the required linear systems.
This makes the implementation of FEAST relatively straight forward, and it ensures that the convergence
rate of FEAST
depends only on the dimension of the subspace being used and on the number of terms $n_c$ in the integration
quadrature rule.

There are many applications of considerable importance, however, where we would like to solve an eigenvalue
problem
by using FEAST, but the use of a direct solver for solving the linear systems is either inadvisable or
impossible.
A direct solver requires that one be able to form and store a factorization of the matrices $(z_kI-A)$.
A recent Parallel FEAST (PFEAST) implementation was proposed for solving larger system sizes of this kind, taking advantage of distributed-memory sparse linear system solvers and domain decomposition techniques \cite{tang2014new,kalantzis2016domain}.
\ep{In very large-scale applications}, however, the structure of the matrix $A$ causes the factorization step to
be extremely slow and expensive to perform,
and the storage of the factorization may even be impossible due memory constraints.
In some other cases, the matrix $A$ is too large and dense to be stored at all, and is instead being represented
implicitly by a rule for performing fast matrix-vector products (see \cite{giannozzi2009quantum} for an example of an application where this approach is used). 
In situations like these, an obvious alternative might be to use iterative linear system solvers rather than direct ones.
With iterative solvers, \bg{assuming that a preconditioner is not used}, one only needs a rule for matrix-vector
multiplication in order to solve a linear system, and there is no need to form or store large, expensive factorizations.

In the following sections, we consider the effectiveness of using iterative linear system solvers when implementing
the FEAST algorithm. In particular, we investigate whether or not the FEAST algorithm can converge quickly and
reliably when the linear systems in the quadrature rule (\ref{eqn:contour_int_quad}) are deliberately solved
inaccurately with
considerable error.  
We also consider the relationship between the resulting modified FEAST algorithm and traditional Krylov subspace
methods for solving eigenvalue problems.




\subsection{Prior Work: Inexact Shift-Invert Subspace Iterations}

\ep{Various} authors \cite{robbe2009inexact,golub2000inexact,berns2006inexact,lai1997inexact} have previously
examined \ep{ the efficiency
  of inner-outer iterations  for solving the eigenvalue problem  using inexact linear system solves for the shift-invert
  subspace iteration procedure}.
Shift-invert subspace iterations find the eigenvectors of a matrix whose eigenvalues are near some shift $\sigma$.

This is done by using subspace iterations with the matrix $(\sigma I-A)^{-1}$, which is `equivalent' to using the
FEAST algorithm with a single
shifted linear system. With {\em inexact} shift-invert subspace iterations, the matrix multiplications $Y=(\sigma I-A)^{-1}X$ are
calculated by solving for $Y$ inexactly using an iterative linear system solver.

The authors in \cite{robbe2009inexact} show that, for general, non-Hermitian matrices, inexact shift-invert subspace
iterations
converge linearly to the eigenpairs of interest, provided that the shifted linear systems are solved sufficiently
accurately
(Theorem 3.1 in \cite{robbe2009inexact}). The required accuracy for the linear systems is an upper bound on the
linear
system residuals that is proportional to the current residual of the eigenvectors; as the eigenvalue problem
converges,
the linear systems must be solved \ep{increasingly} more accurately to ensure convergence. 

They also show that, when using GMRES as the linear system solver,
the number of GMRES iterations that is required to meet the condition for convergence is approximately
the same at each subspace iteration (Proposition 3.8 in \cite{robbe2009inexact}). In other words, although the
shifted linear systems
must be solved to increasing levels of accuracy as the eigenvalue problem converges,
the amount of computation that is required to solve these systems at each subspace iteration generally does not
increase.
This is true without the use of a preconditioner; in fact, most \ep{standard} preconditioning strategies will
prevent this effect from occurring,
thereby increasing the cumulative amount of work that must be done to solve the linear systems.
\ep{We note that the authors in \cite{robbe2009inexact} address this issue by using a tuned preconditioner,
  but we do not consider
  that approach in this paper.}
All of this suggests that using approximate, non-preconditioned linear system solves with shift-invert
subspace iterations
can be a very efficient way to solve an eigenvalue problem.

Like traditional shift-invert iterations, FEAST allows one to find eigenvalues anywhere in the complex plane.\bg{Unlike traditional shift-invert subspace iterations, the convergence rate of FEAST can be systematically improved by changing the number and location of the shifts, and }
\ep{the conditioning of the FEAST shifted matrices can be significantly better because the
  complex shifts can be located farther away
  from the eigenvalues of interest (more particularly if the eigenvalues are located in the real axis).}
By solving its associated linear systems inexactly, we intend to maintain the benefits of using FEAST while taking advantage
of the useful properties of inexact shift-invert subspace iterations.


\section{Iterative FEAST}

``Iterative FEAST" (or IFEAST) is the FEAST algorithm implemented such that the linear systems are deliberately
solved inaccurately. That is, the linear systems are solved \ep{such that the resulting residuals satisfy
a convergence criteria} that is greater (possibly substantially greater) than machine precision.
\ep{Algorithm \ref{algo:ifeast} summarizes the iterative IFEAST procedure for solving the standard Hermitian
  interior eigenvalue problem.}

\begin{algorithm}
\ep{
  \medskip
\noindent \begin{minipage}{\linewidth}
{\footnotesize

\noindent {\bf Start with:}
\begin{itemize}
\item Matrix $A \in \mathbb{C}^{n \times n}$
\item Interval ${\cal I}=(\lambda _{min},\lambda _{max} )$ wherein
fewer than $m_0$ eigenvalues are expected to be found, and closed contour $\cal{C}$ that encloses $\cal I$ in the complex plane
\item Initial guess $X^{(0)} \in \mathbb{C}^{n\times m_0}$ for the search subspace spanned by the solution to the eigenvalue problem
\item Set of $n_c$ quadrature weights and points $(\omega_k,z_k)$ for numerically integrating equation (\ref{eqn:contour_int})
\item Initial value for FEAST eigenvector residual $||R_F||$ \footnote{This can be calculated exactly,
  but we find a good initial value to simply be $||R_F||$=1}
\item Relative tolerance $\alpha$ for linear system residuals, with $0<\alpha<1$
\end{itemize}
\medskip
\noindent {\bf For each} subspace iteration $i$:
\medskip
\begin{description}
\setlength{\itemsep}{3pt}
\item[1.] Iteratively solve $n_c$ shifted linear systems for $Y_k^{(i)}\in \mathbb{C}^{n\times m_0}$.
\begin{equation}\nonumber
(z_kI-A)Y_k^{(i)}=X^{(i)}, \ \ 1\leq k\leq n_c
\end{equation}
\noindent such that the iterations are stopped when the following criterion on the linear system residuals is met:
\begin{equation}\nonumber
||X^{(i)}-(z_kI-A)Y_k^{(i)}||\leq \alpha||R_F||
\end{equation}
\item[2.] Form the filtered subspace $Q$
\begin{equation}\nonumber
Q=\hat \rho(A)X^{(i)} = \sum_{k=1}^{n_c}\omega_k Y_k^{(i)}
\end{equation}

\item[3.] Perform Rayleigh-Ritz procedure to find a new estimate for eigenvalues and eigenvectors:
\begin{description}
\setlength{\itemsep}{3pt}
\item[i.] Solve reduced eigenvalue problem for $X_Q  \in \mathbb{C}^{m_0\times m_0} $
\begin{equation}\nonumber
  A_QX_Q= B_Q X_Q \Lambda
\end{equation}
\noindent with $A_Q =Q^T A Q \text{ and } B_Q=Q^TQ$
\item[ii.] Get new estimate for subspace $X^{(i+1)}$: $X^{(i+1)} = Q X_Q$
\end{description}
\item[4.] Calculate the FEAST eigenvector residual $||R_F||=max\  ||Ax_j-\lambda _j x_j||,\ 1 \leq j \leq m_0,\ \lambda_j \in {\cal I}$.
  If $||R_F||$ is above a given tolerance, {\bf GOTO 1}.
\end{description}
}
\end{minipage}
\caption{\label{algo:ifeast} The IFEAST Hermitian algorithm for solving $AX_I=X_I\Lambda_I$}
}
\end{algorithm}

The implementation of IFEAST requires a new parameter, $\alpha$, that determines the stopping criterion that
is used in solving the linear systems iteratively. Importantly, the stopping criterion changes at each
iteration in proportion
to the eigenvector residuals. IFEAST will not necessarily converge for all values
of $\alpha$ (an issue that we deal with
quantitatively in Section \ref{sec:convergence}), and so it should be heuristically \ep{underestimated}.
We find that, for example, a value of $\alpha=10^{-2}$ tends to work very well in many cases. 
Any iterative linear system solving algorithm that can be used with general matrices can also be used for
solving the linear systems of IFEAST. In this work \ep{we consider only Hermitian matrices $A$}, and
we choose to work with MINRES \cite{paige1975solution} because of its combination
of speed, robustness, and limited storage requirements. \bg{Although the FEAST linear systems $(z_kI-A)$ are not Hermitian, they can still be solved with MINRES because they are shifted versions of Hermitian systems\cite{freund1990conjugate}}.

In the following subsections we describe the properties of IFEAST analytically.  
\ep{We show that IFEAST can converge linearly when its linear systems
are solved inexactly,  and we show how the accuracy of the inexact solves interacts with the other parameters that govern
the behavior of IFEAST.} We also examine the relationship between IFEAST and traditional Krylov eigenvalue solving algorithms;
because the vast majority of \ep{the computation in IFEAST consists of performing calculations with Krylov
  subspaces in order to solve linear systems, it is natural to ask whether or not IFEAST itself is some kind
  Krylov subspace method as a result}.


\subsection{Convergence}
\label{sec:convergence}
The convergence proof from \cite{robbe2009inexact} (Theorem 3.1) can not be straightforwardly applied to \ep{IFEAST}
because IFEAST uses a linear combination of shifted systems, rather than a single shifted system.
We instead offer an alternative proof of convergence \ep{that can be applied} 
to IFEAST for any number of shifts. We show that IFEAST converges linearly by providing an upper bound on
the eigenvector error at a given FEAST subspace iteration that depends linearly on an upper bound on the
eigenvector error at the previous
FEAST subspace iteration. \ep{Below}, we describe this upper bound and the implications that it has for the behavior
and convergence of IFEAST. These results are a modification of the analysis of conventional subspace iterations found
in \cite{saad1992numerical}, and we provide the details of the derivation in \ep{Appendix} \ref{apA}.

For IFEAST, an upper bound on the eigenvector error is given by
\begin{equation}
\label{eqn:upperbound}
||\tilde w_j||\leq \left(\frac{|\gamma_{m_0+1}|+\alpha_j\Delta}{|\gamma_j|}\right)\ ||w_j||.
\end{equation}

The norm $||w_j||$ is an upper bound on the error for the estimation of the eigenvector $x_j$ in the current FEAST
subspace $Q$. The norm $||\tilde w_j||$ is an upper bound on the error for estimating $x_j$ in the FEAST subspace
at the next iteration, $\hat \rho(A)Q$.\bg{The value $\gamma _j$ is the $j^\text{th}$ largest eigenvalue of $\hat\rho(A)$, with corresponding eigenvector $x_j$.}
\ep{The dimension of the FEAST
  search subspace is $m_0$.}

\ep{If all the shifted linear systems are solved inaccurately with a given convergence
criteria $\epsilon$ on the residual norm, i.e.
  ($e_j$ being the unit vector)
\begin{equation}
\label{eqn:feastlinsystol}
||Xe_j-\frac{1}{\omega_k}(z_kI-A)Y_ke_j||\leq\epsilon, \ \forall k,j
 \end{equation}
then the scalar $\alpha_j$ can be defined as the ratio of the magnitude of the maximum FEAST linear system residual $\epsilon$ to the value of $||w_j||$:
\begin{equation}\label{eq:alpha}
\alpha_j=\epsilon/||w_j||.
\end{equation}}
\ep{Also derived in Appendix \ref{apA}}, the scalar $\Delta$ is a function of the spectrum of the matrix that we are diagonalizing,
the locations of the FEAST linear system shifts, and the values of the FEAST linear system weights:
\begin{equation}
\label{eqn:deltadef}
\Delta=\sum_{k=1}^{n_c} ||\omega_k(z_kI-A)^{-1}||.
\end{equation}

The criterion for IFEAST to converge for the eigenvector $x_j$ comes straightforwardly from (\ref{eqn:upperbound}):
\begin{equation}
\label{eqn:convergecondition}
\alpha_j\Delta<|\gamma_j|-|\gamma_{m_0+1}|.
\end{equation}
\ep{Provided that the inequality in (\ref{eqn:convergecondition}) is true, the FEAST subspace will
  become a better approximation
  to the eigenvector subspace of interest with each subsequent iteration. IFEAST will converge
 at the rate of
  $(|\gamma_{m_0+1}|+\alpha_j\Delta)/|\gamma_j|$.
  The smaller the magnitude of this coefficient is, the faster FEAST converges by subspace iteration.
  When the linear systems of IFEAST are solved exactly (i.e. $\alpha_j=0$ at machine precision),
   the convergence rate of traditional FEAST is recovered \cite{tangFEAST}.}
\ep{In turn, if  $\alpha_j$ has the same
  value at every IFEAST iteration and
 (\ref{eqn:convergecondition})  is satisfied, then
  the upper bound (\ref{eqn:upperbound}) guarantees linear convergence.}

For values of $\alpha_j\Delta$ that are much smaller than $|\gamma_{m_0+1}|$, IFEAST behaves similarly
to traditional FEAST:
solving additional linear systems in parallel leads directly to a better convergence rate. If $\alpha_j\Delta$
is on the order of, or greater than,
$|\gamma_{m_0+1}|$, then the behavior of IFEAST \ep{is different from} that of traditional FEAST. In this case,
solving additional linear systems
in parallel does not make IFEAST converge faster, and the convergence rate is dominated by the accuracy of
\ep{ the linear system solves}.
\ep{In addition, we note that the closer the shifts $z_k$ are to the eigenvalues of $A$, the larger $\Delta$
  becomes.
  As a result, the linear systems of IFEAST must be solved to a certain level of accuracy in order to ensure
  that all additional shifted linear systems
  can effectively contribute to a faster convergence rate  (which can be challenging because using more shifted systems means that more complex shifts end up closer to an eigenvalue located on the real axis). Some eigenvalue problems are  also expected to be inherently
more difficult for IFEAST to solve than others, 
\bg{such as,} for example, when there is a cluster of eigenvalues located just outside
the eigenvalue search interval ${\cal I}=(\lambda_\text{min},\lambda_\text{max})$, which can have the effect of
causing the difference
$|\gamma_j|-|\gamma_{m_0+1}|$ from Equation \ref{eqn:convergecondition} to be very small (if $m_0$ is not large
enough).}

\ep{Finally, we point out that
  the definition for $\alpha_j$ in (\ref{eq:alpha}) differs from the definition of $\alpha$ that is used for
  all eigenvector $x_j$
in the IFEAST Algorithm \ref{algo:ifeast}}. This is by necessity,
 as it is not possible to know the maximum
norm of the linear system residuals
that will ensure convergence of the eigenvalue problem \ep{(i.e. condition (\ref{eqn:convergecondition}))}, without
having already solved the eigenvalue
problem. In Algorithm  \ref{algo:ifeast}, we use the value of the FEAST eigenvector residual $||R_F||$ \bg{in place of the eigenvector error $||w_j||$} and
  heuristically determined the parameter $\alpha$ in order to provide
estimates for the linear system residual tolerance.

\subsection{Solving Inexact FEAST Linear Systems with GMRES}
\label{sec:useGMRES}

Unlike the case for convergence, the results from Proposition 3.8 in \cite{robbe2009inexact} can be applied to IFEAST
without modification. That is, in IFEAST, as with inexact shift-invert iterations, the number of GMRES iterations that are
required to satisfy the tolerance on the linear system residuals \ep{(Step 1
in the IFEAST Algorithm \ref{algo:ifeast})} generally does not increase as the eigenvalue problem converges,
even though the tolerance itself becomes smaller at each subsequent IFEAST subspace iteration.
We refer the reader to \cite{robbe2009inexact} for the details of the proof, but point out here that the fundamental reason
for this is simple. The closer the right hand side of a linear system of equations is to an invariant subspace of
the coefficient matrix, the fewer GMRES iterations are required to solve the system to a given tolerance.
As IFEAST iterations converge, the right hand sides of the IFEAST linear systems become closer to being invariant
subspaces
of the matrix that is being diagonalized (and hence become easier to solve), at the same time as the tolerance
for the solution
is made more difficult to reach.

This is also the reason that most linear system preconditioners will actually make the eigenvalue problem
more expensive to solve.
The right hand sides of the IFEAST linear systems converge to invariant subspaces of the matrix that is being
diagonalized, but they
generally do not converge to invariant subspaces of the {\em preconditioned} matrix. In order to use a preconditioner
with IFEAST
without increasing the amount of work that needs to be done to solve the eigenvalue problem, it is necessary to choose
a preconditioner
that either shares the eigenvectors of the matrix being diagonalized, or to choose a new preconditioner at each subspace
iteration such
that the right hand sides of the linear systems are invariant subspaces of the preconditioned matrix. One such strategy
is described
in \cite{freitag2008tuned}.

This effect also makes it difficult to efficiently apply \ep{Algorithm \ref{algo:ifeast}} to generalized eigenvalue
problems (\ref{eqn:eigenvalueproblem}), where the FEAST linear systems become \ep{(at iteration $i$):
\begin{equation}
\label{eqn:generalized_linsys}
(z_kB-A)Y_k^{(i)}=BX^{(i)}.
\end{equation}}
In this case, unlike the standard eigenvalue problem case, the right hand sides do not converge to an invariant
subspace of $(z_kB-A)$, and so the number of GMRES iterations that is required for convergence increases with each
subspace iteration.
It is always possible to rewrite equation (\ref{eqn:generalized_linsys}) so that the right hand sides do converge
  to invariant subspaces
of the coefficient matrix, if we consider solving, for example, $(z_k-B^{-1}A)Y_k^{(i)}=X^{(i)}$.
However, in doing so, we replace our original problem with another problem of at least equal difficulty\bg{; when we replace $(z_kB-A)$ with $(z_k-B^{-1}A)$, every matrix multiplication by $A$ must be accompanied by a linear system solve with $B$, which dramatically increases the cost of iteratively solving the corresponding linear system $(z_k-B^{-1}A)Y_k^{(i)}=X^{(i)}$.}

Several authors have suggested some ways of addressing this challenge. \ep{One possibility involves using tuned
  preconditioners
to recover the desired behavior of GMRES \cite{freitag2008rayleigh,xue2011fast}. Another consists of changing the way
the initial guess is chosen for GMRES \cite{golub2000inexact,ye2011inexact}. }


\section {Relationship between IFEAST and Krylov methods}
\label{sec:feastkrylov}

Standard Krylov eigenvalue solving methods (such as Lanzcos and Arnoldi)
work by building a basis $V$ for the Krylov subspace ${\cal K}(A,X^{(0)})$, using some
initial guess  $X^{(0)}$ for the eigenvectors i.e.

\begin{equation}
\label{eqn:vinkbasis}
V\in {\cal K}(A,X^{(0)})=span\{ X^{(0)},AX^{(0)},A^2X^{(0)},...,A^{k-1}X^{(0)} \},
\end{equation}
\ep{with}
\begin{equation}
X^{(0)}\in \mathbb{C}^{n\times m_0},\ \ V\in \mathbb{C}^{n\times m_0k}.
\end{equation}

For easy comparison with IFEAST below, we consider the case of a block Krylov method where the block size is
\ep{$m_0$ (i.e. size of the FEAST search subspace)}. 
Traditional Krylov methods then use the Rayleigh-Ritz procedure to form and solve a reduced-dimension eigenvalue problem in order to
find approximate eigenpairs in the subspace ${\cal K}(A,X^{(0)})$

\begin{equation}
\label{eqn:reduced_krylov}
\ep{(V^HAV)X_V=(V^HV)X_V\Lambda.}
\end{equation}
\ep{Let us assume that the degree $(k-1)$ of the Krylov subspace (\ref{eqn:vinkbasis})
is made as large as is practically possible.}
\ep{If the residuals of the approximate eigenpairs from the reduced problem (\ref{eqn:reduced_krylov})  do not converge,
then the method can be ``restarted" 
by using a block of Ritz vectors $X^{(1)}$ from the solution of (\ref{eqn:reduced_krylov})
as the starting vectors for building a new Krylov subspace ${\cal K}(A,X^{(1)})$ of degree $(k-1)$.}

FEAST,  when doing the contour integration exactly, forms a subspace by applying a spectral projector
to $X^{(0)}$, which is then also used to solve a reduced-dimension eigenvalue problem i.e.

\begin{equation}
\label{eqn:contint_reg}
Q=\rho (A)X^{(0)}= \frac{1}{2\pi i}\oint_{\cal C} (zI-A)^{-1}X^{(0)}dz,
\end{equation}

\begin{equation}
\label{eqn:feastreducedevp}
\ep{(Q^HAQ)X_Q=(Q^HQ)X_Q\Lambda.}
\end{equation}

We can understand the relationship between FEAST and traditional Krylov methods by considering what happens
when the integrand $(zI-A)^{-1}X^{(0)}$ in (\ref{eqn:contint_reg}) is evaluated approximately by using a
Krylov subspace. We can rewrite the integral (\ref{eqn:contint_reg}) as:

\begin{equation}
\label{eqn:qfuncy}
Q=\rho(A)X^{(0)}=\frac{1}{2\pi i}\oint_{\cal C} Y(z)dz,
\end{equation}

\noindent where $Y(z)$ is the solution to the linear system

\begin{equation}
\label{eqn:yzlinsys}
(zI-A)Y(z)=X^{(0)}.
\end{equation}

\noindent If we use a Krylov subspace method to find an approximate solution to (\ref{eqn:yzlinsys}), then

\begin{equation}\label{eqn:yz}
  \ep{Y(z)=VY_V(z),\ \ Y_V(z)\in \mathbb{C}^{m_0k\times m_0}},
\end{equation}

\noindent where $V$ is the same Krylov subspace basis from equation (\ref{eqn:vinkbasis}),
and $Y_V$ \ep{is an approximate solution to
$(zI-A)VY_V(z)=X^{(0)}$.}
Importantly, the Krylov basis $V$ is not a function of $z$, because the
Krylov subspace that is generated by $(zI-A)$ depends only on the matrix $A$ and not on the shift $z$.
Because $V$ is independent of $z$, we can rewrite the expression (\ref{eqn:qfuncy})
for $Q$ in such a way that the FEAST reduced-dimension eigenvalue
problem (\ref{eqn:feastreducedevp}) takes a familiar form. Rewriting the expression for $Q$, we get

\begin{equation}
\label{eqn:contint_krylov}
Q= \frac{1}{2\pi i}\oint_{\cal C} Y(z)dz = VG_v,
\end{equation}
with
\begin{equation}
\label{eqn:gexp}
G_v\in \mathbb{C}^{m_0k\times m_0}=\frac{1}{2\pi i}\oint_{\cal C} Y_V(z)dz.
\end{equation}
\noindent Then, the FEAST reduced eigenvalue problem (\ref{eqn:feastreducedevp}) becomes

\begin{equation}
\label{eqn:feastreducedevpkrylov}
(G_v^HV^HAVG_v)X_Q=(G_v^HV^HVG_v)X_Q\Lambda.
\end{equation}

Comparing (\ref{eqn:feastreducedevpkrylov}) with (\ref{eqn:reduced_krylov}) makes it clear that IFEAST
itself is, in fact, a Krylov subspace method. The difference between IFEAST and more traditional Krylov
methods is that IFEAST \bg{uses contour integration to} select an \ep{ideally-suited} linear combination of vectors
from the Krylov basis $V$ for
finding the desired eigenvalues, without first having to solve a reduced eigenvalue problem in that basis.

\bg{Being able to select the desired eigenvalues in this way can have substantial benefits. One of the challenges in using Krylov subspaces is that finding certain eigenvalues, particularly interior eigenvalues or eigenvalues that are clustered closely together, can require a subspace basis $V$ of \ep{very large} dimension. Using a large-dimension subspace basis $V$ entails large storage requirements for that basis, and a large computational cost for solving the corresponding reduced eigenvalue problem (\ref{eqn:reduced_krylov}). When using IFEAST, on the other hand, the dimension of the reduced eigenvalue problem (\ref{eqn:feastreducedevpkrylov}) is always $m_0$, which is substantially smaller than the dimension $km_0$ of the traditional reduced eigenvalue problem (\ref{eqn:reduced_krylov}).

Moreover, when IFEAST is implemented with a linear system solver that uses a short recurrence relation (e.g. MINRES), then it can solve eigenvalue problems by using a Krylov subspace of arbitrarily large dimension without having to form and store a basis for that subspace; by using short recurrences, IFEAST can form the $n\times m_0$ matrix product $Q=VG_v$ without forming or storing either the $n\times km_0$ matrix $V$ or the $km_0\times m$ matrix $G_v$.} Thus, eigenpairs that would previously have been difficult or impossible to \ep{obtain} due to constraints on the dimension of $V$ become much more tractable to calculate,
and  the spectrum slicing capability of FEAST is maintained by making it possible to selectively \ep{find specific eigenpairs anywhere in the spectrum.}




\bg{The relationship between IFEAST and traditional Krylov methods also offers a different perspective on achieving convergence when using restarts. In the context of IFEAST, a Krylov restart amounts to an approximate subspace iteration with $\rho (A)$ for a particular choice of contour $\cal C$. Using contour integration to choose the subspace with which to restart ensures that restarting will reliably result in convergence, with inequality (\ref{eqn:upperbound}) giving quantitative answers regarding whether or not restarting will result in convergence and, if it does, how quickly convergence will occur. IFEAST reverses the process that is used in other restarting strategies \cite{sorensen1992implicit,saad1984chebyshev}, in which the subspace that is used for restarting is determined {\em after} solving a reduced eigenvalue problem in the full Krylov subspace, rather than before.}

\bg{We elaborate further on the relationship between IFEAST and traditional Krylov techniques in the following subsections, where we show how the implementation of IFEAST with particular linear system solvers is related to other Krylov subspace methods for solving eigenvalue problems. We show that implementing IFEAST using the Full Orthogonalization Method (FOM) is equivalent to traditional explicitly restarted block Arnoldi, and that implementing IFEAST using GMRES is closely related to using Harmonic Rayleigh-Ritz for interior eigenvalue problems. }

\bg{




\subsection{IFEAST + FOM is Restarted Arnoldi}
The block Arnoldi method constructs an orthonormal basis $V$ $\in \mathbb{C}^{n\times m_0k}$ of block size $m_0$ and Krylov polynomial degree $k-1$ (for a total dimension of $m_0k$), and then solves a reduced eigenvalue problem from the Rayleigh-Ritz method in order to find estimates for the desired eigenvalues and eigenvectors, i.e.

\begin{equation}
HX_V=X_V\Lambda
\end{equation}

\begin{equation}
H=V^HAV, \ \ \ V={\cal K}_{k}(A,X^{(0)}), \ \ V^HV=I
\end{equation}

\noindent where $H\in \mathbb{C}^{m_0k\times m_0k}$ is upper Hessenberg and $X^{(0)}\in \mathbb{C}^{n\times m_0}$ is the initial guess for the eigenvectors. If the residuals on the estimated eigenpairs $(VX_V,\Lambda)$ are not good enough, then the method can be explicitly ``restarted'' by building a new Krylov subspace ${\cal K}_{k}(A,X^{(1)})$ using a new starting block $X^{(1)}$. The new starting block consists of linear combinations of the estimated eigenvectors, i.e.

\begin{equation}
  \label{eqn:restartvecs}
X^{(1)}=VX_VM
\end{equation}

\noindent where $M\in \mathbb{C}^{m_0k\times m_0}$ gives the linear combinations that are used to determine each vector in the new starting block. A variety of different choices for $M$ are possible \cite{saad1992numerical}. A single iteration of IFEAST, when implemented with FOM, produces a new estimate for the eigenvectors of interest $X^{(1)}$ that is equivalent to expression (\ref{eqn:restartvecs}) for a particular, natural choice of $M$.

Implementing IFEAST requires forming a subspace $Q\in \mathbb{C}^{n\times m_0}$ by evaluating the contour integral (\ref{eqn:qfuncy}), which in turn requires solving linear systems of the form (\ref{eqn:yzlinsys}). We restate these tasks (respectively) here, i.e.

\begin{equation}
\label{eqn:qfuncy2}
Q=\rho(A)X^{(0)}=\frac{1}{2\pi i}\oint_{\cal C} Y(z)dz,
\end{equation}

\begin{equation}
\label{eqn:yzlinsys2}
(zI-A)Y(z)=X^{(0)}.
\end{equation}

\noindent FOM is used to solve the linear system (\ref{eqn:yzlinsys2}) by forming $V$ using Arnoldi iterations, and then solving a projected linear system \cite{saad2003iterative}, i.e.

\begin{equation}
Y(z)=V \left( V^H(zI-A)V \right)^{-1}V^HX^{(0)}.
\end{equation}

\noindent Because the linear system matrix $(zI-A)$ is just a shifted version of the original matrix $A$, the solution for $Y(z)$  can be written in terms of the upper Hessenberg matrix that is generated by the Arnoldi method, i.e.

\begin{equation}
  \label{eqn:yzfom}
Y(z)=V (zI-H)^{-1}V^HX^{(0)}.
\end{equation}

\noindent Inserting this into the expression for the IFEAST subspace $Q$ (\ref{eqn:qfuncy2}) , it becomes clear that using FOM is equivalent to applying the FEAST filter function $\rho(\lambda)$ to the upper Hessenberg matrix $H$ from Arnoldi

\begin{equation}
\label{eqn:qfuncH}
Q=V\frac{1}{2\pi i}\oint_{\cal C}(zI-H)^{-1}dz V^HX^{(0)}=V\rho(H)V^HX^{(0)}.
\end{equation}

This is equivalent to filtering out the components of the unwanted Arnoldi Ritz vectors from $X^{(0)}$, leaving only the Ritz vectors whose Ritz values are inside the contour $\cal C$ in the complex plane. We can see this by writing the eigenvalue decomposition of $H$ and reordering its eigenvalues and eigenvectors so that the wanted eigenpairs (i.e. the ones whose eigenvalues are inside $\cal C$) are grouped together, i.e.

\begin{equation}
  H=X_V\Lambda X^H_V,
\end{equation}

\begin{equation}
X_V=\left[ X_w\ \ X_u \right], \ \ \ \Lambda=\left[\begin{matrix}\Lambda_w & 0\\ 0 & \Lambda_u \end{matrix}\right],
\end{equation}

\noindent and by writing the initial guess $X^{(0)}$ in terms of its Ritz vector components in the $V$ subspace

\begin{equation}
X^{(0)}=VX_wW+VX_uU,
\end{equation}

\noindent where $(X_w,\Lambda_w)$ are the $m_0$ wanted Ritz eigenpairs (i.e. the ones whose eigenvalues are inside $\cal C$ in the complex plane), $(X_u,\Lambda_u)$ are the $(k-1)m_0$ unwanted Ritz eigenpairs, and $W$ and $U$ are the components of $X^{(0)}$ in terms of the wanted and unwanted Ritz eigenvectors (respectively). Rewriting (\ref{eqn:qfuncH}) in these terms, we get

\begin{align}
\label{eqn:qfuncH2}
Q&=V \left[ X_w\ \ X_u \right] \left[\begin{matrix}\rho(\Lambda_w) & 0\\ 0 & \rho(\Lambda_u) \end{matrix}\right] \left[ X_w\ \ X_u \right]^H  V^H (VX_wW+VX_uU) \\
&=V(X_w\rho(\Lambda _w)W+X_u\rho(\Lambda _u)U).
\end{align}

IFEAST with FOM thus forms a subspace by filtering the Ritz values and vectors from the Arnoldi Rayleigh Ritz matrix $H$; the components of $X^{(0)}$ in the direction of the wanted Ritz vectors are kept roughly the same, and the components of $X^{(0)}$ in the direction of the unwanted Ritz vectors are substantially reduced. When the contour integral in (\ref{eqn:qfuncH}) is evaluated exactly, then $\rho(\Lambda _w)=I_{m_0\times m_0}$ and $\rho(\Lambda _u)=0_{(k-1)m_0\times (k-1)m_0}$, and IFEAST forms and solves a reduced eigenvalue problem using only the Arnoldi Ritz vectors corresponding to the wanted Ritz values. The vectors that are used as the initial guess for the next IFEAST iteration, then, are just the normalized Arnoldi Ritz vectors corresponding to the Ritz values that are inside the contour $\cal C$ in the complex plane, i.e.

\begin{equation}
  \label{eqn:ifeastrestart}
X^{(1)}=VX_w=VX_V\left[ \begin{matrix} I_{m_0\times m_0} \\ 0_{(k-1)m_0\times m_0} \end{matrix}\right].
\end{equation}

\noindent IFEAST with FOM is equivalent, then, to performing block Arnoldi with a restart strategy that consists of selecting the desired Ritz vectors and discarding the rest.

In practice this restart strategy can be unreliable for obtaining eigenvalues in the interior of the spectrum. One perspective on why this happens is that the Rayleigh-Ritz procedure works well for resolving exterior eigenvalues, but not for resolving interior ones; restarting with Ritz vectors is thus unreliable for obtaining interior eigenvalues \cite{morgan1998harmonic}. A remedy for this is to use the Harmonic Rayleigh Ritz procedure \cite{morgan1991computing,morgan1998harmonic}, wherein one solves a different reduced eigenvalue problem that more accurately obtains the eigenvalues that are located near some shift.

The fact that the restart strategy (\ref{eqn:ifeastrestart}) is equivalent to using FOM with IFEAST suggests another perspective on why it is ineffective. Getting  IFEAST to converge requires solving its associated linear systems such that their residuals are sufficiently small, and FOM does not minimize the linear system residual for a given subspace. Reliably achieving convergence for interior eigenpairs requires the use of a linear system solver that minimizes the linear system residual, such as GMRES or MINRES.

\subsection{IFEAST + GMRES is related to Harmonic Rayleigh Ritz}

In fact, using GMRES with IFEAST is closely related to using the Harmonic Rayleigh Ritz procedure.
When using GMRES to solve (\ref{eqn:yzlinsys2}) for $Y(z)$, the solution takes the form \cite{saad2003iterative}

\begin{equation}
Y(z)=V \left( V^H(zI-A)^H (zI-A)V \right)^{-1}V^H(zI-A)^HX^{(0)},
\end{equation}

\noindent where $V$, again, is the block Arnoldi basis. The IFEAST subspace $Q$ then becomes

\begin{equation}
  \label{eqn:gmresq}
Q=V\left( \frac{1}{2\pi i} \oint_{\cal C} \left[ V^H(zI-A)^H(zI-A)V\right]^{-1}V^H(zI-A)^HV  dz  \right)  X^{(0)}_V,
\end{equation}

\noindent where $VX^{(0)}_V=X^{(0)}$ is the initial guess $X^{(0)}$ expressed in the Arnoldi basis $V$.

The integrand in (\ref{eqn:gmresq}) is equivalent to the matrix that one arrives at when using Harmonic Rayleigh Ritz with Arnoldi. With Harmonic Rayleigh Ritz, one seeks to find approximations for the eigenvalues that are near some shift $z\in \mathbb{C}$, using the subspace basis $V$. This is done by solving the reduced, generalized eigenvalue problem \cite{morgan1991computing,morgan1998harmonic}

\begin{equation}
  \label{eqn:genharmonicritz}
A_V(z)X_V(z)=B_V(z)X_V(z)(zI-\Lambda (z)),
\end{equation}

\begin{equation}
A_V(z) = V^H(zI-A)^H(zI-A)V, \ \ \ \ B_V(z)=V^H(zI-A)^HV,
\end{equation}

\noindent where $VX_V(z)$ are now the Harmonic Ritz vectors, and $\Lambda (z)$ are the Harmonic Ritz values. In most applications the shift $z$ is taken to be a fixed parameter, but here we are considering a case where it will vary, making the projected matrices $A_V(z)$ and $B_V(z)$, and the Harmonic Ritz vectors and values $X_V(z)$ and $\Lambda (z)$, into matrix-valued functions of the shift. Like any generalized eigenvalue problem, (\ref{eqn:genharmonicritz}) can be written as a standard, non-symmetric eigenvalue problem with a corresponding eigenvalue decomposition, i.e.

\begin{equation}
\label{eqn:harmonic_decomp}
B^{-1}_V(z)A_V(z)=X_V(z)(zI-\Lambda (z))X_V^{-1}(z).
\end{equation}


\noindent If we note that

\begin{equation}
\left[B_V^{-1}(z)A_V(z)\right]^{-1}=\left[ V^H(zI-A)^H(zI-A)V\right] ^{-1}V^H(zI-A)^HV,
\end{equation}

\noindent then we can use this combined with Equation (\ref{eqn:harmonic_decomp}) in order to write the expression for $Q$ (\ref{eqn:gmresq}) in terms of the Harmonic Rayleigh Ritz eigenvalue decomposition:

\begin{align}
  \label{eqn:gmresqharm}
Q&=V\left( \frac{1}{2\pi i} \oint_{\cal C} \left[ B^{-1}_V(z)A_V(z) \right]^{-1} dz  \right)  X^{(0)}_V,\\
&=V\left( \frac{1}{2\pi i} \oint_{\cal C} \left[ zI-X_V(z)\Lambda (z) X_V^{-1}(z) \right]^{-1} dz  \right)  X^{(0)}_V. \label{eqn:rhononlin}
\end{align}

Generating the IFEAST subspace by using GMRES is thus equivalent to using contour integration to filter the initial guess by using Arnoldi Harmonic Ritz values and vectors. Unlike with FOM, however, the resulting contour integral is not equivalent to applying the usual FEAST spectral filter $\rho(\lambda)$ to a projected matrix. Instead, the integration in (\ref{eqn:rhononlin}) is the contour integral of the resolvent of a {\em nonlinear} eigenvalue problem, where the eigenvalues and eigenvectors are functions of the complex variable $z$ that are derived from the Harmonic Rayleigh Ritz procedure.

}


\section{Results and \ep{Discussions}}
\label{section:results}


In this section we illustrate the behavior of IFEAST using two example matrices.

\subsection{Example I: Si$_2$}
\begin{figure}
  \begin{center}
  \textbf{Integration Contour and Linear Shifts for Si$_2$}
\end{center}
  \par\medskip
\includegraphics[width=0.98\textwidth]{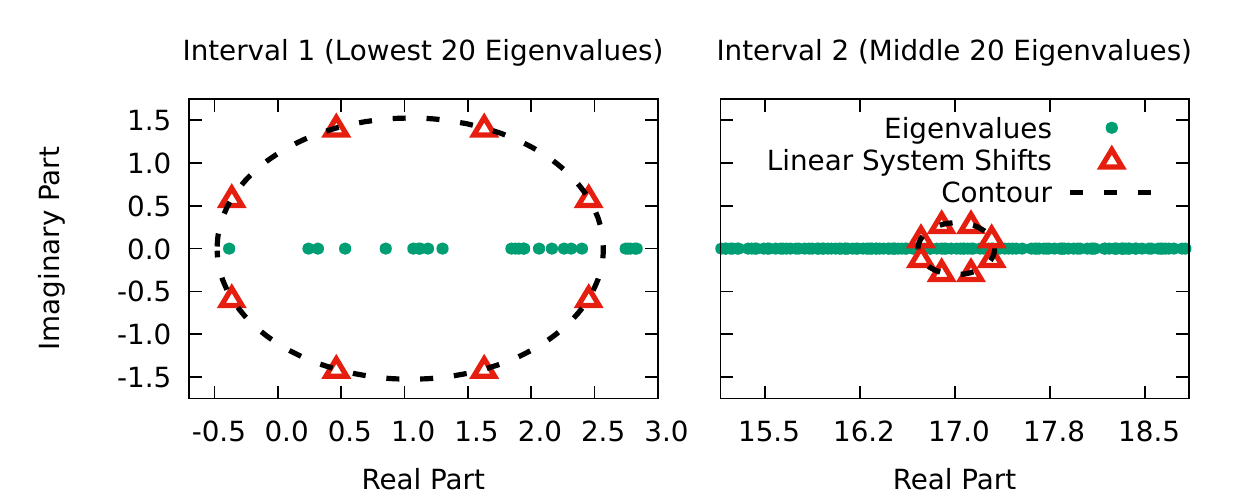}
\caption{\label{fig:contours} Locations in the complex plane of the IFEAST integration contour and linear system shifts.
  Contour and shifts are provided for two calculations: finding the lowest 20 eigenvalues, and finding
  the middle 20 eigenvalues.
  \ep{The quadrature nodes (shifts) are located on a perfect circle
    (although the contour appears elliptical due to the bounds of the plots).
  We use here a total of 4 linear system shifts for discretizing the integral in the upper-half contours using
  the Trapezoidal rule.}
}
\end{figure}

Our first example is the Si$_2$ matrix from the University of Florida Sparse
Matrix Collection \cite{timdaviscollection}. Si$_2$ is a \ep{real symmetric} $769\times 769$
matrix from the electronic structure code PARSEC; it represents the Hamiltonian operator
of a quantum system consisting of two silicon atoms.
We illustrate the behavior of IFEAST by
calculating eigenvector/eigenvalue pairs in two places in the spectrum of Si$_2$: the
lowest 20 eigenpairs, and the middle 20 eigenpairs.
The eigenvalues, search
contours, and linear system shifts for each of these calculations are illustrated in Figure~\ref{fig:contours}.
\ep{Using the same scale, we note that the contour for Interval 1 (the lowest eigenvalues)
is much larger than the contour for Interval 2 (the middle eigenvalues) because the eigenvalues in Interval 2
are clustered much more closely together.
   Due to the symmetry property of FEAST for addressing the Hermitian problem \cite{polizziFEAST},
  it is only necessary to perform the numerical quadrature on
  the upper-half of the contour by using $n_{c_{up}}=n_c/2$ total shifted linear systems.
The trapezoidal rule is used to select the location and weight of each of these shifts.}

\ep{For a given number of quadrature points $n_c$,  subspace size $m_0\geq m$ ($m$ the number of eigenvalues, here 20), 
equation (\ref{eqn:upperbound}) tells us that a small enough convergence criterion for the linear systems $\alpha$ guarantees that the FEAST linear convergence criteria depends entirely
  on the the value of the filtering function i.e. the outer-iteration subspace iteration.
For the Si$_2$ example, in particular, if we select $m_0=1.5m=30$ and $\alpha\leq10^{-1}$, IFEAST converges in 9
outer-iterations for both contours (using $n_{c_{up}}=4$).
In the following examples we deliberately choose parameter values such that the behavior of IFEAST deviates from that of conventional FEAST, in order to illustrate the effects of inexact linear sytem solves.}

\begin{figure}
  \begin{center}
  \textbf{Si$_2$: Subspace Iteration Convergence}
\end{center}
  \par\medskip
\includegraphics[width=0.98\textwidth]{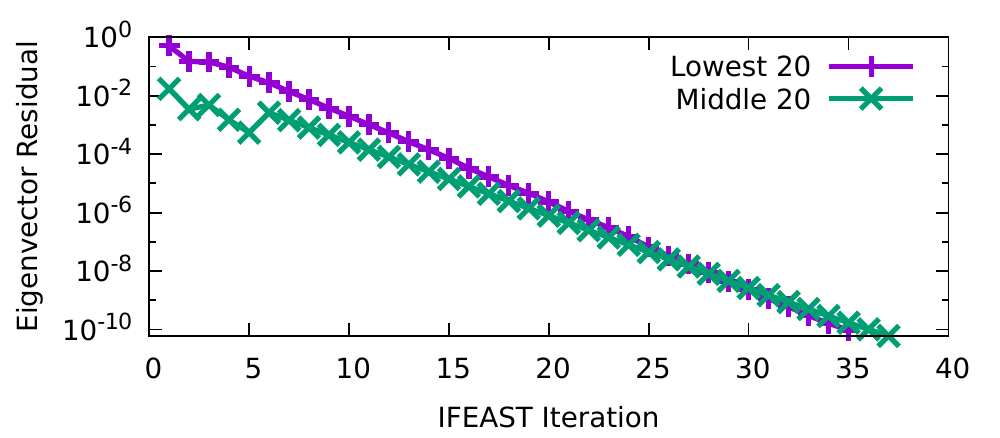}
\caption{Convergence of IFEAST eigenvalue calculations for Si$_2$ by subspace iteration. MINRES is used as the linear system solver,
  with a subspace dimension of $m_0=20$ and a rather large linear system convergence criterion of $\alpha=1/2$.}
\label{fig:si2convsi}
\end{figure}
\begin{figure}
  \begin{center}
  \textbf{Si$_2$: MINRES Iterations at each IFEAST Iteration}
\end{center}
  \par\medskip
\includegraphics[width=0.98\textwidth]{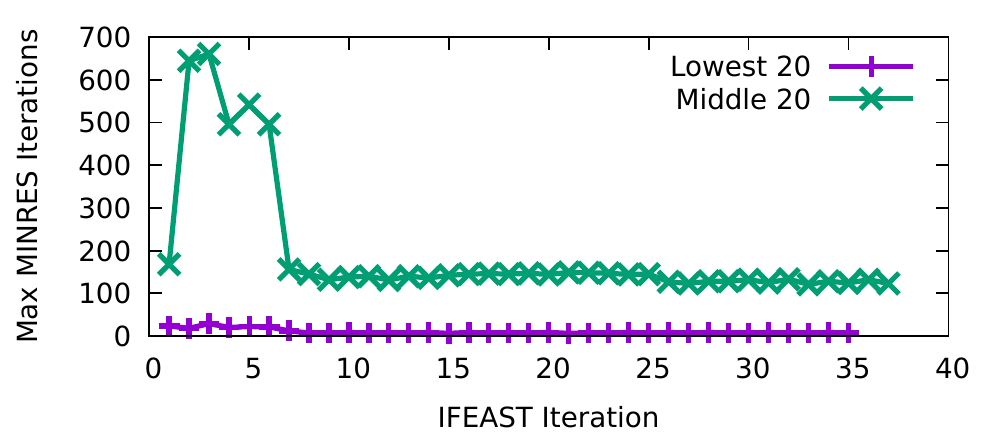}
\caption{Maximum number of MINRES iterations performed at each IFEAST iteration for Si$_2$, for calculating both the lowest 20
  eigenvalues and the middle 20 eigenvalues (see Figure \ref{fig:contours}). The maximum number of MINRES iterations is the number of MINRES
  iterations required by the shifted linear system right hand side that takes the longest to converge. When using enough parallelism to solve
  all right hand sides simultaneously, this is a measure of the amount of time that each IFEAST iteration takes.}
\label{fig:si2convmatvecIt}
\end{figure}

\ep{
  Figure \ref{fig:si2convsi} shows the eigenvector residual at each
subspace iteration when using the IFEAST Algorithm with MINRES as the linear system solver, the smallest possible subspace
size of $m_0=20$, and the linear system convergence criterion $\alpha=1/2$. Interval 1 and Interval 2
converge at similar rates by subspace iteration and, as expected, the convergence for each is linear.}
Unlike with traditional FEAST, however, the number of subspace iterations that is required for convergence
is not a good measure of the amount of time \ep{needed by IFEAST for solving the} eigenvalue problem.
\ep{Indeed}, when solving the linear systems iteratively, some shifted systems will converge faster than others,
and some right hand sides will converge faster than others. If enough parallel processing power is available to solve all linear
system right hand sides simultaneously, then the best measure of the amount of time that a single IFEAST iteration takes is
the number MINRES iterations that is required for the most difficult linear system right hand side to converge. This is shown in
Figure \ref{fig:si2convmatvecIt}. As \ep{specified}  by Proposition 3.8 in \cite{robbe2009inexact}, the number of MINRES iterations required at
each FEAST iteration is approximately constant.

\begin{figure}
  \begin{center}
  \textbf{Si$_2$: Matrix Vector Product Convergence}
\end{center}
  \par\medskip
\includegraphics[width=0.98\textwidth]{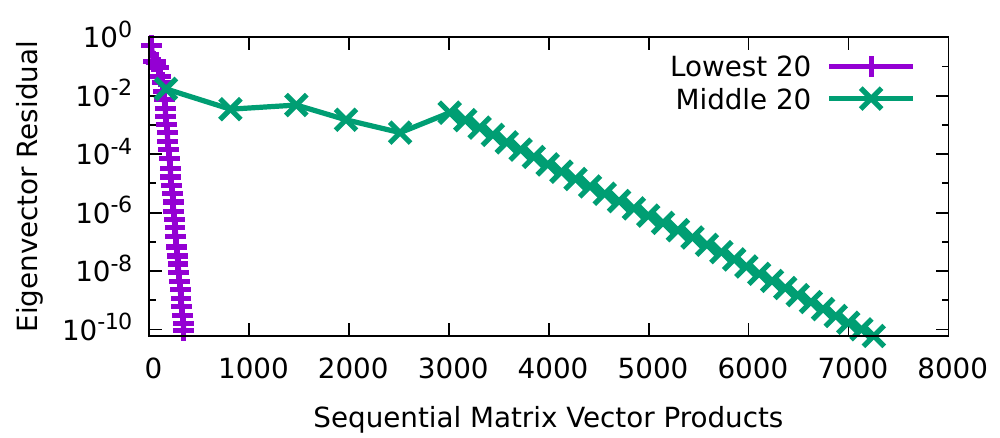}
\caption{Convergence of IFEAST eigenvalue calculations for Si$_2$ by number of sequential matrix vector multiplications, for both the lowest
  20 eigenvalues and the middle 20 eigenvalues (see Figure \ref{fig:contours}). The number of sequential matrix vector multiplications
  is the sum of the number of MINRES iterations
  for the slowest-converging linear system right hand side at each IFEAST iteration. This is the best measure of the amount of time
  that IFEAST requires to converge when solving all linear system right hand sides in parallel at each subspace iteration.}
\label{fig:si2convmatvec}
\end{figure}

Figure \ref{fig:si2convmatvec} displays this information in a different way, showing the cumulative number of sequential matrix
vector products that is required to reach a given eigenvector residual when all linear system right hand sides are solved
in parallel with MINRES. The time required for convergence of the eigenvalue problem is proportional to the number
of sequential matrix vector products, and so Figure \ref{fig:si2convmatvec} gives the best comparison of the
performance of IFEAST for Interval 1 and Interval 2. 

\ep{It is clear from the results in Figure \ref{fig:si2convmatvec} that
  the convergence for Interval 1 happens much quicker than for Interval 2.}
The reason for this is that the eigenvalues both inside and around Interval 2 are closely clustered together in the middle of the spectrum,
whereas the eigenvalues in Interval 1 are well-separated at the lower edge of the spectrum.
Any Krylov subspace algorithm will require many more iterations to find the eigenvalues in Interval 2 than
it will for Interval 1, and IFEAST is no exception. However, one benefit of implementing IFEAST with MINRES
is that, unlike with other Krylov eigenvalue methods, the size of the subspace needed for calculating the eigenpairs
in Interval 2 is exactly the same as the size of the subspace needed for calculating the eigenpairs in Interval 1.
This makes it possible to maintain the parallelism of traditional FEAST by solving for many eigenpairs in parallel by using multiple contours.

\begin{figure}
  \begin{center}
  \textbf{Si$_2$: Subspace Iteration Convergence for \\ Different Numbers of Shifted Systems}
\end{center}
  \par\medskip
\includegraphics[width=0.98\textwidth]{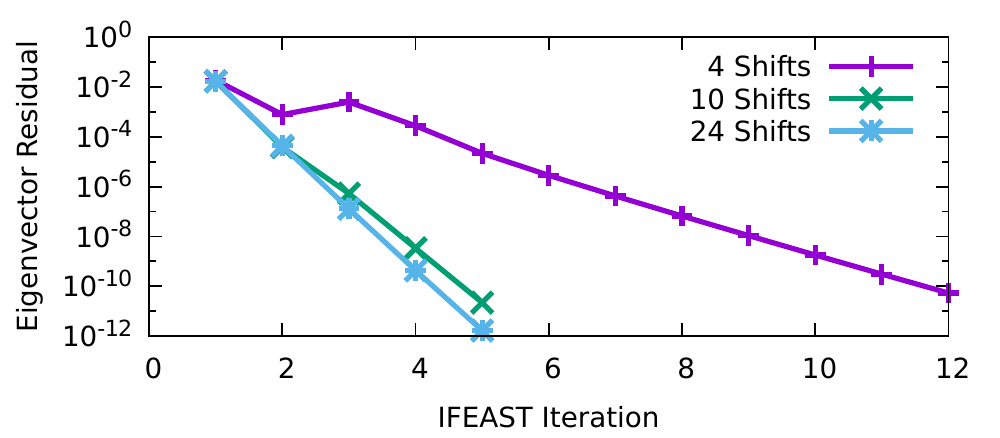}
\caption{ Eigenvector residual versus IFEAST iteration for calculating the 20 eigenvalues in Interval 2
  (the middle eigenvalues of Si$_2$; see Figure \ref{fig:contours}), using several different numbers of shifted linear systems $n_{c_{up}}$
  in the upper-half contour, and for a constant subspace size of $m_0=25$. Increasing the number
  of shifted systems that are solved in parallel improves the convergence rate up to
  the point that convergence becomes limited by the
  accuracy of the linear system solves.}
\label{fig:si2convsicp}
\end{figure}

In traditional FEAST, the rate of convergence by subspace iteration can always be improved by increasing the accuracy
of the numerical integration of Equation (\ref{eqn:contour_int}), usually by increasing the number of terms in the quadrature
rule (\ref{eqn:contour_int_quad}). As discussed in Section \ref{sec:convergence}, the situation is less simple in IFEAST,
due to the relationship between convergence and the accuracy of the linear system solutions. In general, increasing the
number of shifted linear systems in the quadrature rule (\ref{eqn:contour_int_quad}) will improve the convergence rate
by subspace iteration up to the point that convergence becomes limited by the accuracy of the linear system solutions,
after which increasing the number of quadrature points will no longer improve convergence.
\ep{This effect is illustrated in Figure \ref{fig:si2convsicp}, which shows the convergence of IFEAST
by subspace iteration for several numbers of linear system shifts}. \ep{Here we calculate the 20 eigenpairs inside Interval 2
(see Figure \ref{fig:contours}) by using a subspace size of $m_0=25$, with MINRES again as the linear system solver and a
convergence criterion $\alpha=1/2$. With these parameters, increasing the number of shifted linear systems $n_{c_{up}}$ from 4 to 10
increases the subspace iteration convergence rate considerably, but increasing the number of shifted linear systems from
10 to 24 barely changes the convergence rate at all.}

\begin{figure}
  \begin{center}
  \textbf{Si$_2$: Matrix Vector Product Convergence for \\ Different Numbers of Shifted Systems}
\end{center}
  \par\medskip
\includegraphics[width=0.98\textwidth]{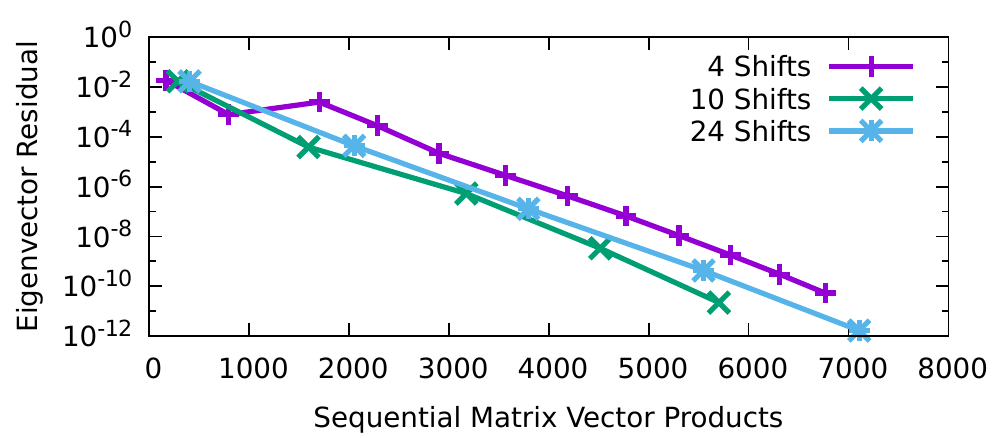}
\caption{Eigenvector residual versus number of required sequential matrix vector multiplications for calculating
  the 20 eigenvalues in Interval 2 (the middle eigenvalues of Si$_2$; see Figure \ref{fig:contours}), using several
  different numbers of shifted linear systems. All shifted linear system right hand sides are assumed to be solved in parallel.
  Increasing the number of shifted linear systems can cause IFEAST to take longer to converge when doing so brings some of those
  shifts closer to the eigenvalues of the matrix without also increasing the subspace iteration convergence rate.}
\label{fig:si2convmatveccp}
\end{figure}

As \ep{previously} mentioned, the convergence rate by subspace iteration is not a true measure of the performance of IFEAST;
a better measure of performance is the number of matrix vector products that is required to reach a given eigenvector residual when
all linear system right hand sides are solved in parallel. Figure \ref{fig:si2convmatveccp} shows the convergence of IFEAST versus
the number of sequential matrix vector products for the same numerical experiment, where we calculate the 20 eigenpairs
inside Interval 2 by using several different numbers of shifted linear systems. When looking at the required number of matrix vector products,
increasing the number of shifted linear systems from 4 to 10 improves performance, but increasing the number of shifted linear systems
from 10 to 24 actually decreases performance, resulting in IFEAST taking longer to converge.
\ep{Although the  shifted linear systems are solved in parallel, some of the linear systems
  are more difficult to solve than others, because their shifts are closer to the real axis (and are thus closer to the eigenvalues).
  Due to the limited
  accuracy of the linear system solves, the convergence rate by subspace iteration remains essentially
  the same for both 10 and 24 shifts.
As a result, IFEAST requires more time when using 24 shifts because it needs to do more work to solve the linear systems that are closer to the real axis, while at the same time the limited accuracy of the linear system solves prevents it from converging more quickly.}

\subsection{Example II: Na$_5$}
Our second example is the Na$_5$ matrix, also from the University of Florida Sparse Matrix Collection \cite{timdaviscollection}.
Na$_5$ is a real symmetric $5832\times 5832$ matrix from the electronic structure code PARSEC, and it represents the Hamiltonian operator
of a quantum system consisting of five sodium atoms. In order to provide context for the performance of IFEAST, we calculate the 50
lowest eigenvalue/eigenvector pairs for Na$_5$ using both IFEAST and Arnoldi. The implementation of Arnoldi that we use is ARPACK
\cite{arpacksite}, which implements single vector Arnoldi with implicit restarts.

\begin{table}
  \begin{center}

    \textbf{Na$_5$: Matrix Vector Product Comparison for IFEAST and Arnoldi}

    \par\medskip
  \begin{tabular}{ | l | c | c | c | c |}
    \hline
    \textbf{Subspace Size} & 75 & 100 & 200 & 787 \\ \hline \hline

    \textbf{IFEAST Iterations}  & 11 & 11 & 11 & 10 \\ \hline
    \textbf{Arnoldi Restarts} & 57 & 23 & 5 & 0 \\ \hline \hline

    \textbf{IFEAST Total Matvec}  & 34,819 & 44,212 & 105,283 & 463,691 \\ \hline
    \textbf{IFEAST Sequential Matvec} & 672 & 498 & 368 & 278 \\ \hline
    \textbf{Arnoldi Matvec} & 946 & 854 & 844 & 787 \\ \hline

  \end{tabular}
\end{center}
\medskip
\caption{Comparison of the number of matrix vector products required to calculate the 50 lowest eigenpairs of Na$_5$ to an eigenvector
  residual of $10^{-10}$, using both IFEAST and Arnoldi (from the package ARPACK). Matrix vector product counts are shown for several subspace
  sizes. The IFEAST ``Total Matvec" is the total number of matrix vector products that IFEAST requires, and the IFEAST ``Sequential Matvec"
  is the number of matrix vector products that must be done sequentially if all of the matrix vector products that can be done in parallel
  are performed in parallel. }
  \label{tab:arvsfeast}
\end{table}

Table \ref{tab:arvsfeast} shows the number of matrix vector products that is required to calculate the 50 lowest eigenvalue/eigenvector pairs
of Na$_5$ to an eigenvector residual of $10^{-10}$ using both IFEAST and Arnoldi. The results are shown for several different values of subspace size.
For IFEAST, the subspace size is the value of the parameter $m_0$ in the IFEAST algorithm, which is the size of the subspace that is used for the FEAST
subspace iterations. For Arnoldi, the subspace size is the maximum size of the Krylov subspace for ARPACK.
IFEAST is run by using 4 linear system shifts in the upper-half
contour, that are chosen by using the trapezoidal rule, and a linear system
convergence criterion of $\alpha=1/5$. MINRES is used as the linear system solver.

IFEAST generally requires a substantially larger total number of matrix vector products than Arnoldi does.
However, if the available parallelism in the IFEAST algorithm is fully utilized, \ep{meaning that all matrix vector products 
are performed in parallel}, 
then the relative number of matrix vector products that must be done sequentially \ep{becomes competitive in comparison with Arnoldi}. 
In other words, although IFEAST has to do much more work than Arnoldi, it can be made to converge faster in time by doing most of that work in parallel.

\begin{figure}
  \begin{center}
  \textbf{Na$_5$: MINRES Iterations at each IFEAST Iteration}
\end{center}
  \par\medskip
\includegraphics[width=0.98\textwidth]{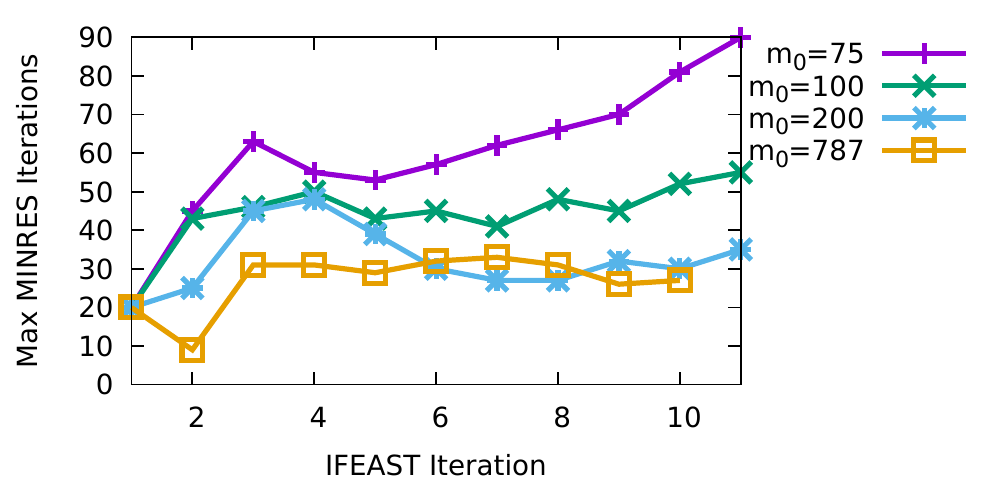}
\caption{Maximum number of MINRES iterations over all linear system right hand sides at each IFEAST iteration when calculating
  the lowest 50 eigenpairs of Na$_5$, using several different subspace sizes. The maximum number of MINRES iterations
  over all right hand sides is a measure of the amount of time that each IFEAST iteration
  requires when all linear system right hand sides are solved in parallel. Using larger subspace sizes reduces the number
  of required MINRES iterations by allowing IFEAST to select better Ritz pairs to use as approximations to the desired eigenpairs;
  this is equivalent to choosing better Ritz pairs with which to restart a Krylov subspace algorithm.
 }
\label{fig:na5minresit}
\end{figure}

\ep{Table \ref{tab:arvsfeast} also compares the number of IFEAST iterations
to the number of Arnoldi restarts;
as described in Section \ref{sec:feastkrylov}, IFEAST subspace iterations are equivalent to restarts when a Krylov algorithm is used for
solving the linear systems. We note that the number of Arnoldi restarts is greater than the number
of IFEAST iterations for smaller subspace size (i.e.  $m_0=75$ and $m_0=100$),
however, increasing the Krylov subspace size reduces the number of restarts
that is required for Arnoldi. The number of iterations (i.e. restarts) for IFEAST, on the other hand, is roughly the same regardless of the IFEAST subspace size.
The performance of IFEAST keeps improving with increasing subspace size because IFEAST is then able
to select better Ritz vectors with which to restart, thereby reducing the required number of linear system iterations (i.e. the degree of
the Krylov polynomial) at each subspace iteration.} This effect is illustrated in Figure \ref{fig:na5minresit}, which shows the maximum number
of MINRES iterations over all linear system right hand sides at each IFEAST iteration when calculating the lowest 50 eigenpairs of Na$_5$ using
several different subspace sizes.

\section{Conclusion}

By implementing the FEAST eigenvalue algorithm using iterative linear system solvers, it is possible to quickly and robustly
calculate select eigenpairs of a matrix anywhere in its spectrum by using only matrix vector multiplication. The resulting eigenvalue
algorithm, which we call Iterative FEAST (IFEAST), is equivalent to a block Krylov subspace algorithm that uses contour integration
in order to determine the linear combinations of Krylov basis vectors that are used for restarting. 

Depending on the particular linear system solving algorithm that is used in its implementation, IFEAST can be shown to be mathematically equivalent to other, well-known Krylov eigenvalue algorithms. The distinguishing feature of the IFEAST algorithm is that, in its actual implementation, it uses Krylov subspaces indirectly for solving linear systems of equations, rather than directly for projecting the original eigenvalue problem. IFEAST can thus take advantage of linear system solving algorithms like MINRES or BICGSTAB that do not need to store a basis for the Krylov subspace. This makes it possible to solve an eigenvalue problem by implicitly using a very large Krylov subspace without ever having to store a basis for it. 


In being able to solve for select eigenvalue/eigenvector pairs using an almost arbitrarily small amount of storage, IFEAST retains the spectrum slicing property of traditional FEAST, making it possible to solve for large numbers of eigenpairs in parallel. The IFEAST algorithm also retains many of the other parallel characteristics of traditional FEAST (such as solving multiple linear systems, and multiple right hand sides, in parallel),
with the caveat that the benefits of this parallelism can be diminished when the shifted linear systems are solved too inaccurately relative
to the accuracy of quadrature rule for the FEAST contour integration.

As described in Section \ref{sec:useGMRES},
future research directions will include applying the work of previous authors \cite{golub2000inexact,freitag2008rayleigh,xue2011fast,ye2011inexact} to IFEAST in order to try to make it as efficient for the generalized eigenvalue problem as it is for the standard eigenvalue problem.




\appendix

\section{FEAST Convergence Bounds}
\label{apA}

In this section we show how to derive the upper bound on the eigenvector error for inexact FEAST

\ep{
\begin{equation}
\label{eqn:upperbound2}
||x_j-\tilde q_j||\leq \left(\frac{|\gamma_{m_0+1}|+\alpha_j\Delta}{|\gamma_j|}\right)\ ||x_j- q_j||.
\end{equation}
}

The upper bound in (\ref{eqn:upperbound2}) can be derived using a modification of the method
for finding an upper bound on the eigenvector error that is used in analyzing standard subspace iterations \cite{saad1992numerical}.

\subsection{Standard Subspace Iterations}

With standard subspace iterations, we want to find the eigenvectors corresponding
to the $m_0$ largest-magnitude eigenvalues of a matrix $A\ \in \mathbb{C}^{n\times n}$. This is done
by repeatedly multiplying an approximate subspace $Q\in \mathbb{C}^{n\times m_0}$ by $A$, and reorthogonalizing
the column vectors of $Q$ in between multiplications \ep{(using Rayleigh-Ritz, for example)}. The usual method for proving
convergence is to show that, for every eigenvector $x_j,\ 1\leq j\leq m_0$, an upper bound on the error of
its estimation in the subspace $Q$ goes down after each subspace iteration. This can be done by judiciously
\ep{choosing} a vector $q_j \in Q$ that is close to $x_j$ and showing that there is always a different vector
$\tilde q_j \in AQ$ that is closer to $x_j$ than $q_j$ is.

Let $X_1 \in \mathbb{C}^{n\times m_0}$ be the subspace whose column vectors are the eigenvectors that we want
to find, and $X_2 \in \mathbb{C}^{n\times (n-m_0)}$ be the subspace composed of the other $n-m_0$ eigenvectors.
Then the vector $q_j$ is usually chosen to be the unique vector in $Q$ that satisfies

\begin{equation}
X_1X_1^Tq_j=x_j.
\end{equation}

\noindent In that case, the difference vector $w_j=q_j-x_j$ is spanned exactly by $X_2$,
since $X_1$ and $X_2$ are mutually orthogonal, invariant subspaces of A. The vector $\tilde q_j$ is then chosen to be

\begin{equation}
\label{eqn:newsivec}
\tilde q_j=\frac{1}{\lambda_j}Aq_j,
\end{equation}

\noindent where $\lambda_j$ is the eigenvalue corresponding to the eigenvector $x_j$. The difference vector
$\tilde w_j=\tilde q_j-x_j$ is then also spanned exactly by $X_2$, a fact that we can
use to relate $||w_j||$ to $||\tilde w_j||$ i.e.

\begin{align}
&\tilde q_j=\frac{1}{\lambda _j}Aq_j = \frac{1}{\lambda _j}(Ax_j+Aw_j)=x_j+\frac{1}{\lambda_j}Aw_j,\\
&\tilde w_j=\tilde q_j-x_j=\frac{1}{\lambda_j}Aw_j,\\
&||\tilde w_j||=\frac{1}{|\lambda_j|}||Aw_j||\leq \frac{|\lambda_{m_0+1}|}{|\lambda_j|}||w_j||, \label{eqn:si_upperbound}
\end{align}
where we know that $||Aw_j||\leq|\lambda_{m_0+1}|||w_j||$ because $w_j$ is spanned exactly by $X_2$,
the $(n-m_0)$ eigenvectors corresponding to the eigenvalues with magnitudes less than or equal
to $|\lambda_{m_0+1}|$.
Equation (\ref{eqn:si_upperbound}) shows that an upper bound on the error for the estimation
of $x_j$ in the subspace $Q$ always decreases when $Q$ is multiplied by $A$, and that it
does so at a rate that is linear and proportional to the ratio between $\lambda_{m_0+1}$ and
$\lambda_j$. Thus, subspace iterations are guaranteed to converge faster when the subspace
$Q$ is larger and when the eigenvalues of $A$ are more separated.

\subsection{Inexact FEAST}

We can find a similarly informative upper bound with which to analyze the convergence of
iterative FEAST by following a similar line of reasoning.
Traditional FEAST can be interpreted as a subspace iteration that uses the matrix $\hat \rho (A)$
instead of the original matrix $A$,

\begin{equation}
A \longrightarrow \hat\rho(A)=\sum_{k=1}^{n_c} \omega _k (z_kI-A)^{-1}.
\end{equation}

\noindent Then the upper bound (\ref{eqn:si_upperbound}) becomes

\begin{equation}
||\tilde w_j||\leq \frac{|\gamma_{m_0+1}|}{|\gamma _j|}||w_j||, \label{eqn:feast_upperbound}
\end{equation}
\bg{where $\gamma _j$ is the $j^\text{th}$ largest eigenvalue of $\hat\rho(A)$, with corresponding eigenvector $x_j$.}
  \ep{The $\gamma_j$  with the largest  magnitudes correspond to the `wanted' eigenvalues of $A$ that lie inside of the integration
  contour (\ref{eqn:contour_int}).} Making the quadrature rule (\ref{eqn:contour_int_quad})
more accurate by increasing the number of quadrature points $n_c$ has the effect of making the ratio
$|\gamma_{m_0+1}|/|\gamma_{j}|$ smaller, which is how standard FEAST can improve its rate of
convergence by solving more linear systems.

Equation (\ref{eqn:feast_upperbound}) requires modification when the linear systems of FEAST are solved inexactly.
In particular, if we apply $\hat\rho(A)$ by solving the linear systems
\ep{
\begin{equation}
\label{eqn:linsys}
\frac{1}{\omega_k}(z_kI-A)y_{k,j}=q_j, \ \forall k=1,\dots, n_c, \ \forall j=1,\dots,m_0
\end{equation}

\noindent such that there is some error $s_{k,j}$ in the solution of the linear system, i.e.} 

\begin{equation}
s_{k,j}=\omega_k(z_kI-A)^{-1}q_j-y_{k,j},
\end{equation}

\noindent then, for inexact FEAST, equation (\ref{eqn:newsivec}) becomes 
\ep{
\begin{equation}
\label{eqn:newfeastvec}
\tilde q_j=\frac{1}{\gamma_j}\left(\hat\rho(A)q_j-\sum_{k=1}^{n_c}s_{k,j}\right).
\end{equation}
}
\ep{This is not necessarily very useful in practice, however, because the values of $s_{k,j}$ are
not known.} Instead, (\ref{eqn:newfeastvec}) can be rewritten in terms of the linear system residuals,
the norms of which are used as the stopping criteria for iterative linear system solvers:
\ep{
\begin{equation}\label{eqn:qj}
\tilde q_j=\frac{1}{\gamma_j}\left(\hat\rho(A)q_j-\sum_{k=1}^{n_c}\omega_k(z_kI-A)^{-1}r_{k,j}\right),
\end{equation}
with
\begin{equation}
r_{k,j}=q_j-\frac{1}{\omega_k}(z_kI-A)y_{k,j}.
\end{equation}
Since $q_j=w_j+x_j$, we can derive the expression for $\tilde w_j$ from (\ref{eqn:qj}):
\begin{equation}
\tilde w_j=\tilde q_j-x_j=\frac{1}{\gamma_j}\left(\hat\rho(A)w_j-\sum_{k=1}^{n_c}\omega_k(z_kI-A)^{-1}r_{k,j}\right).\\
\end{equation}
We can then find an upper bound similar to (\ref{eqn:feast_upperbound}):
\begin{equation}
||\tilde w_j||\leq\frac{|\gamma_{m_0+1}|}{|\gamma_j|}||w_j||+\frac{1}{|\gamma_{j}|}\sum_{k=1}^{n_c}||\omega_k(z_kI-A)^{-1}|| \ ||r_{k,j}|| .\label{eqn:infeast_upperbound}
\end{equation}
Assuming that all the linear systems (\ref{eqn:linsys}) are solved using iterative solvers with the same convergence
criteria $\epsilon$ on the residual norm, then $||r_{k,j}||\leq \epsilon, \ \forall k,j$, we get: 
\begin{equation}
||\tilde w_j||\leq \left(\frac{|\gamma_{m_0+1}|+\alpha_j \Delta}{|\gamma_j|} \right)||w_j||,
\end{equation}
with
\begin{equation}
\alpha_j=\epsilon/||w_j||,
\end{equation}
and
\begin{equation}
\Delta=\sum_{k=1}^{n_c}||\omega_k(z_kI-A)^{-1}||.
\end{equation}
}
\bg{If the linear systems are solved such that $\alpha_j$ is the same at every FEAST subspace iteration,
  then linear convergence is guaranteed, with the rate of convergence depending on accuracy of the linear system solutions.}


  \bigskip\noindent\textbf{Acknowledgment:} We thank  Alessandro Cerioni and Antoine Levitt for bringing
  their FEAST-GMRES
  experiment results to our attention. We thank Yousef Saad and Peter Tang for useful discussions. 
  This material is supported by Intel and NSF under Grant \#CCF-1510010.

\bibliographystyle{siamplain}
\bibliography{mybib}

\end{document}